\documentclass[a4paper,UKenglish,cleveref, autoref, thm-restate]{lipics-v2021}

\usepackage[dvipsnames]{xcolor}
\definecolor{my_green}{HTML}{bbe7c7}
\definecolor{my_orange}{HTML}{f3c589}
\definecolor{my_red}{HTML}{f78e62}
\definecolor{my_pink}{HTML}{f573b8}
\definecolor{my_blue}{HTML}{2585ba}
\definecolor{blue_fig_1}{HTML}{4285f4}
\definecolor{orange_fig_1}{HTML}{f6b26c}

\theoremstyle{definition}

\title{An iterative Constraint Programming approach to integrate maximum workload constraints in preemptive jobshop scheduling} 

\titlerunning{Integrating maximum workload constraints in pJSP} 

\author{Tanguy {Terrien}\footnote{Corresponding author}}{LAAS-CNRS, Université de Toulouse, France  \and \url{https://www.laas.fr/fr/annuaire/756} }{tterrien@laas.fr}{https://orcid.org/0009-0000-6791-1237}{}

\author{Cyrille Briand}{LAAS-CNRS, Université de Toulouse, France \and \url{https://www.laas.fr/fr/annuaire/188} }{briand@laas.fr}{https://orcid.org/0000-0003-1890-9100}{}

\authorrunning{T. Terrien and C. Briand} 

\Copyright{Tanguy Terrien and Cyrille Briand} 



\category{Short paper} 

\relatedversion{} 

\nolinenumbers 

\EventEditors{Nicolas Beldiceanu}
\EventNoEds{1}
\EventLongTitle{32nd International Conference on Principles and Practice of Constraint Programming (CP 2026)}
\EventShortTitle{CP 2026}
\EventAcronym{CP}
\EventYear{2026}
\EventDate{July 20--23, 2026}
\EventLocation{Lisbon, Portugal}
\EventLogo{}
\SeriesVolume{379}
\ArticleNo{60}

\begin{document}

\maketitle

\begin{abstract}

Optimizing schedules in real-world settings often requires considering workload constraints, especially for human resources, to ensure regulatory compliance, impose rest periods, or level the workload over the working horizon. This paper focuses on tackling this family of constraints in the context of preemptive jobshop scheduling, as preemption is particularly relevant when human resources are involved (allowing personnel to flexibly switch between tasks). Preemption also offers theoretical insights as a relaxation of non-preemptive problems. The main contribution of this paper is a Constraint Programming approach designed to handle effectively maximum workload constraints in a preemptive setting, without decomposing activities into unit-duration tasks (which may be computationally prohibitive). Since workload constraints introduce significant additional complexity, we further propose a method that iteratively introduces the workload constraints into the problem, along with tailored heuristics specifically designed to guide the search efficiently. The experimental results demonstrate the effectiveness of our approach on a large set of instances, highlighting its performance compared to a well-known industrial solver, IBM's CP Optimizer{.}\footnote{This work was supported by the French National Research Agency (ANR) under the project \textbf{HIS$^3$} ANR-22-CE10-0012-05. Special thanks to  \href{https://scholar.google.com/citations?user=ph1JJ_IAAAAJ}{\underline{E.Hébrard}} for his help on Mistral.}
\end{abstract}

\vspace{-0.2cm}
\section{Introduction}
\vspace{-0.1cm}

Optimizing schedules in real-world settings requires taking workload constraints into account, especially in personnel scheduling, in order to ensure regulatory compliance, prevent fatigue, and maintain productivity \cite{BerghBBDB13, OzderOE20}. Several sectors such as project management, healthcare (e.g., nurse rostering \cite{AbdelghanyYE21}), and transportation (e.g., railway crew scheduling \cite{HeilHB20}) are subject to such regulations, which makes the integration of workload and rest-time rules a challenging aspect of scheduling. General overviews of scheduling theory and manufacturing processes also emphasize the importance of these considerations \cite{Pinedo2012}. 

Considering an operator (a resource), a maximum workload constraint states that the number of working shifts assigned to the operator must not exceed a given limit over a specified time period. In practice, several such constraints may coexist and interact. For example, an operator may not be allowed to work more than two shifts within any three-shift time window and no more than five shifts during any five-day week. Handling such workload regulations often requires allowing task preemption: an operator can interrupt a task and resume it later. In the scheduling literature, preemption is commonly used to relax otherwise difficult scheduling problems. For example, the preemptive version of the single-machine scheduling problem with release dates and deadlines $[r_j,d_j]$, denoted $1|r_j,pmtn|L_{\max}$, can be solved in polynomial time using the well-known \emph{Jackson Preemptive Schedule} (JPS)~\cite{CarlierPinson1990}. However, this tractability does not extend to more complex scheduling environments. In particular, the preemptive variant of the Jobshop Scheduling Problem (pJSP), which is considered in this paper, remains NP-hard.

Dealing simultaneously with preemption and workload constraints is even more challenging. Common Constraint Programming (CP) or Mixed-Integer Linear Programming (MILP) approaches, such as the ones used in~\cite{Laborie_cpo_2018}, model tasks and rest periods as unit-length intervals. Global constraints are then used to enforce schedule consistency, such as \textbf{NoOverlap} and \textbf{AtMostSeqCard}~\cite{siala_optimal_2014}, which handle disjunctive resources and workload regulation rules, respectively. However, as pointed out in~\cite{Juvin2023}, the resulting expansion of the solution space caused by task decomposition often leads to a degradation in performance, making optimal solutions difficult to obtain in realistic settings.

To avoid the aforementioned combinatorial issue, we introduce in this paper a new class of constraints that capture both preemption and rest-time requirements: the \textbf{MaxW constraints}, standing for \textbf{Max}imum \textbf{W}orkload, and explain how to handle them without decomposing tasks. The main contributions of this paper are as follows:\\
(i) we introduce workload constraints for human operators, generalizing rules such as forbidding consecutive working shifts or enforcing minimum rest periods;\\
(ii) we propose a CP model for the preemptive Job-Shop Scheduling Problem with MaxW constraints (MaxW-pJSP) using the Mistral solver \cite{mistral};\\
(iii) we design an iterative solution approach that progressively activates MaxW constraints when needed, using a lazy constraint generation mechanism guided by dedicated heuristics;\\
(iv) we implement a baseline CP model for the MaxW-pJSP using IBM\textsuperscript{\tiny{\textregistered}} CPO~\cite{Laborie_cpo_2018};\\
(v) we generate new benchmark instances by augmenting classical \href{https://people.brunel.ac.uk/~mastjjb/jeb/orlib/jobshopinfo.html}{\textit{pJSP datasets}} with varying densities of workload constraints, and report extensive computational experiments.

The remainder of the paper is organized as follows. Section~2 reviews the related literature. Section~3 defines the problem, introduces the proposed MaxW constraints and describes our CP model using Mistral as well as the iterative method. Section~4 presents and analyzes the experimental results. Section~5 concludes and outlines directions for future research.

\vspace{-0.2cm}
\section{Brief state-of-the-art}
\label{sota}
\vspace{-0.1cm}

The resource-constrained project scheduling problem (RCPSP) provides a framework to model limited resource availability in scheduling. When human resources are involved, traditional RCPSP formulations are extended to incorporate specific workload constraints such as minimum rest periods between tasks \cite{Laborie_cpo_2018}. MILP is also common for modeling workload constraints, as seen in nurse rostering problems, where integer programming techniques \cite{Santos2016} or column generation-based heuristic approaches~\cite{Strandmark2020} are applied. Exact solutions for large instances can be computationally costly. Metaheuristics are widely employed for intractable large-scale problems such as complex scheduling environments, like nurse rostering \cite{Liu2018, Haspeslagh2014}. In this paper, the focus is on the JSP. It is a classic NP-hard problem that is focused on processing jobs on machines. Introducing human operators with specific constraints makes the problem significantly more complex. Each operation may require not only a machine but also an operator with limited working hours. While the general JSP has been extensively studied, research papers explicitly integrating human resources are scarce. In~\cite{mauguiere_new_2005}, a branch-and-bound algorithm is proposed to deal with operators' unavailability constraints. A so-called filter-and-fan based heuristic is also proposed for scheduling in flexible job shops under workforce constraints in~\cite{Muller03082022}. 
The preemptive Job-Shop Scheduling Problem (pJSP) further complicates the classical JSP by allowing operations to be interrupted and resumed. The pJSP with workload constraints can be seen as a generalization of the already NP-hard pJSP (considering makespan minimization), and is therefore also NP-hard. Introducing workload constraints on human resources intertwines personnel scheduling requirements with the intricate precedence and resource dependencies of job-shop operations. 

Related work on online printing shop scheduling, which also considers flexible job-shop problems with resumable operations and machine unavailability, has relied on CP models \cite{LUNARDI2020105020}. Although this work focuses on makespan minimization without explicitly modeling workload constraints, it provides a useful foundation for addressing pJSP variants with additional operational constraints, as considered in this paper. Sequence-based constraints can also be used to model rest regulations. In particular, the constraint $\textsc{AtMostSeqCard}(X,V,q,u,d)$, where $X = \langle x_1,\ldots,x_n \rangle$ is a sequence of variables and $V$  a set of values, ensures that $\forall i \in \{1,\ldots,n-q+1\}, \quad \sum_{j=i}^{i+q-1} [x_j \in V] \le u$ and $\sum_{j=1}^{n} [x_j \in V] = d$. Such a constraint can, for instance, enforce an operator to have at least $2$ shifts of rest between two working shifts, while still controlling the total number of working shifts over the entire planning horizon.

The above-mentioned works notably assume that tasks or rest periods have unit length. More recently, Juvin et al.~\cite{Juvin2023} proposed an efficient CP approach for the pJSP that avoids this limitation and prevents the combinatorial enumeration of all possible task–shift combinations. In the following, we show how their approach can be extended to handle MaxW constraints.

\vspace{-0.2cm}
\section{Modelling and solution approaches}
\vspace{-0.1cm}

A MaxW constraint associated with an operator $k$ is defined as a pair $(\delta^k,[u,v])$, where $[u,v]$ denotes an interval of shifts and $\delta^k$ is the maximum number of shifts during which operator $k$ can work within this interval. Let $x_t^k$ be a binary variable equal to $1$ if operator $k$ works during shift $t$ and $0$ otherwise. The constraint can then be expressed as $\sum_{t \in [u,v]} x_t^k \leq \delta^k$.

The MaxW-pJSP environment is defined by a set of jobs $\mathcal{J}$. Each job $J_i \in \mathcal{J}$ consists of a sequence of $n_i$ tasks $t_{i,j}$ for $j \in \{1,\dots,n_i\}$. Tasks within a job must satisfy precedence constraints, i.e., $t_{i,j+1}$ can start only after the completion of $t_{i,j}$. Each task $t_{i,j}$ is characterized by a processing time $P_{i,j}$, expressed in number of working shifts, an earliest starting time $s_{i,j}$, and a latest completion time $e_{i,j}$. We consider a set of $K$ operators indexed from $0$ to $K-1$. Each task is assigned to a specific operator who must execute it. Since we consider the preemptive version of the problem, an operator may interrupt a task and resume it later without penalty, provided that the same operator continues its execution. An operator can not process multiple task simultaneously. In addition, each operator $k$ is subject to a set of MaxW constraints. The objective is to minimize the makespan, denoted $c_{\max}$.

To improve computational efficiency, we avoid explicitly modeling every working shift. Instead, we rely on the \textsc{PreemptiveNoOverlap} constraint implemented in the Mistral solver \cite{Juvin2023}. This constraint assigns earliest start and latest end of tasks such that all preemptive tasks assigned to an operator can be fully processed within the available time horizon. Mistral's dedicated propagation mechanisms, including overload checks, allow the pJSP with makespan minimization to be solved efficiently. Once optimal starts and ends of tasks are obtained, an explicit shift-based schedule can be reconstructed as a post-processing step using the classical JPS algorithm \cite{jackson_algo_74}. This approach significantly reduces the number of variables by avoiding modelling every shift individually. We therefore propose a new Mistral-based model for the MaxW-pJSP problem.
\newpage

\textbf{The Mistral model} \quad In the following, if not specified, consider $i\in \{0, ..., |\mathcal{J}|-1\}$, $j\in \{0, ..., n_i\}$, $k \in [0,K-1]$. Uppercase symbols correspond to problem data, and lowercase symbols to decision variables. We refer to MaxW$^k$ as the set of MaxW constraints of operator $k$, and to T$_k$ as the set of tasks of operator $k$ (fixed beforehand). We write $\delta_c^k$ and $[u^k_c,v^k_c]$ the values of MaxW constraint $c$ related to operator $k$.
\vspace{-0.3cm}
$$    \min c_{max}$$
\vspace{-0.9cm}
\begin{align*}
    s_{i,j} \in [0, \mathit{UB} - \text{P}_{i,j}], \text{ } e_{i,j} \in [\text{P}_{i,j}, \mathit{UB}] && \forall i,j \quad (\text{V}1) \\
    c_{max} \in [0, \mathit{UB}]&&  \quad (\text{V}2) \\
    \text{d}^k_q \in [0, |q|]&&\quad \forall k, \forall q \quad (\text{V}3) \\
    e_{i,j} \geq s_{i,j} + \text{P}_{i,j} && \quad \forall i, \forall j \quad (\text{C}1)\\
    s_{i,j+1} \geq e_{i,j} && \quad \forall i, \forall j  \quad (\text{C}2) \\
    s_{i,0} = 0 \quad \text{and} \quad e_{i, n_i} \leq c_{max} && \quad \forall i  \quad (\text{C}3) \\
    \textsc{Pre.NoOverlap}\Bigl(\underset{_{\forall i, \forall j \in \text{T}_k}}{\bigl\{(s_{i,j}, e_{i,j}, \text{P}_{i,j})\bigr\}} \cup
    \underset{\forall q}{\bigl\{(S_{q}, E_{q}, \text{d}^k_q)}\bigr\}\Bigr) && \quad  \forall k \quad (\text{C}4) \\
    \left(v_c^k - u_c^k\right) - \delta_c^k \leq \sum_{ \forall q \in [u^k_c, v^k_c]} \text{d}^k_q && \forall k, \forall c \in \text{MaxW}^k \quad (\text{C}5) \\
    \text{d}^k_q \leq \max_{\substack{\forall c \in \text{MaxW}^k \\ \text{\tiny{st. }} q \in [u^k_c,v^k_c]}} \left(v_c^k - u_c^k\right) - \delta_c^k && \quad \forall k, \forall q \quad (\text{C}6)
\end{align*}

\textbf{Variables:} For each task $j$ of each job $i$, we define a start time variable and an end time variable (V1), where $\mathit{UB}$ is a precomputed upper bound on the schedule length. The makespan $c_{\max}$ is defined in V2. To handle workload constraints for an operator $k$, we partition the scheduling horizon into a set of subintervals $\mathcal{Q}^k$. Every time point at which a MaxW constraint may start or end defines a boundary of one such subinterval. Consequently, the horizon is divided into a finite sequence of non-overlapping time windows that together cover all shifts. For each operator $k$ and each subinterval $q \in \mathcal{Q}^k$, we create an additionnal task with fixed start and end ($S_q$ and $E_q$), and variable duration. Duration $d_q^k$ represents the number of rest shifts that will be assigned to operator $k$ within subinterval $q$. When created, a variable $d_q^k$ ranges from $0$ to $|q|$, where $|q|$ denotes the length (in shifts) of subinterval $q$.

\begin{example}
Let us give an illustrative example. Consider an operator with two MaxW constraints \textcolor{blue_fig_1}{(5, [0, 6])} and \textcolor{orange_fig_1}{(2, [4, 9])}. Figure~\ref{fig:example_MaxW} depicts both MaxW, their computed subintervals $q_1 \text{=} [0, 4],\text{ } q_2 \text{=} [4, 6],\text{ } q_3 \text{=} [6, 9]$, with associated duration variables $d_1, d_2, d_3$ and their corresponding constraints:
\vspace{-0.1cm}
\begin{figure}[htb]
\begin{center}
    \includegraphics[width=13cm, height=3.6cm]{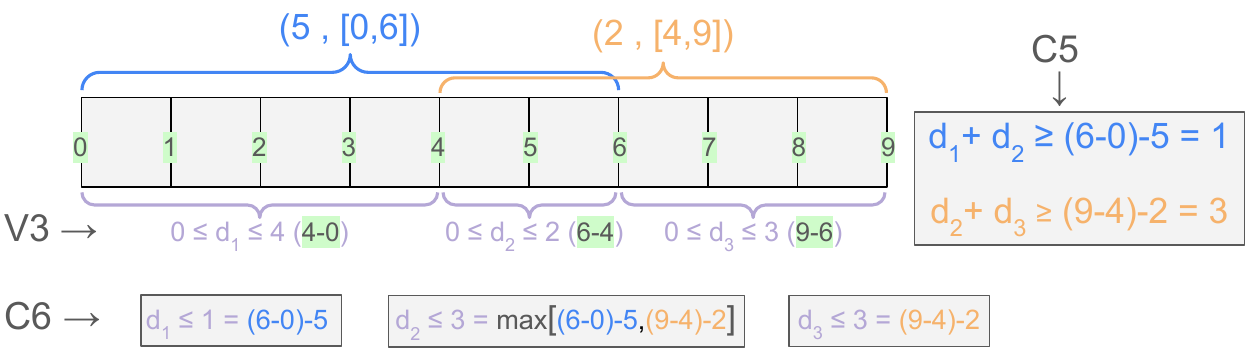}
    \caption{Two MaxW constraints, their subintervals and associated duration variables}
    \label{fig:example_MaxW}
\end{center}
\end{figure}
\end{example}

\textbf{Constraints:} Each task $t_{i,j}$ is subject to a duration constraint ensuring that it is fully executed (C1). Every task is assigned to a specific operator. Tasks belonging to the same job are linked by precedence constraints (C2). Constraints (C3) maintain the consistency of the schedule horizon ($c_{\max}$). To ensure that tasks do not overlap, a \textsc{PreemptiveNoOverlap} constraint is enforced for each operator $k$ (C4), considering both operator $k$'s work tasks and additional rest tasks (with variable durations $d_q^k$). For each operator $k$ and each MaxW constraint, the model identifies the set of rest subintervals contained in the corresponding MaxW interval. The total rest duration assigned to these subintervals must be greater than or equal to the required amount, ensuring compliance with the workload limitation (C5). Finally, constraints (C6) reduce the search space by tightening domains of variables $d_q^k$. Without loss of optimality, $d_q^k$ is upper-bounded by the largest $\delta^k$ among all MaxW constraints containing subinterval $q$. Indeed, assigning this amount of rest shifts is enough to satisfy all corresponding MaxW constraints, and adding more rest shifts could not improve the makespan and is therefore unnecessary. Note that Mistral's original propagation algorithm assumes fixed task durations in the \textsc{PreemptiveNoOverlap} constraint. We therefore adapted it for our model by substituting the minimum duration value of each variable-duration rest task $d_q^k$ in the overload-checking and edge-finding procedures. It may underestimate the mandatory workload of rest variables but never produces incorrect inferences. This extension will be integrated in Mistral's open source code \cite{mistralrepo}. It allows the solver to compute variables $s_{i,j}$, $e_{i,j}$ and $\text{d}^k_q$ in a way that guarantees the existence of a feasible preemptive schedule that is makespan-optimal (produced by the JPS, explained at the start of this section).

\textbf{The iterative approach} \quad To reduce the computational burden induced by the full set of MaxW  constraints, we adopt a lazy constraint generation scheme where the CP model plays the role of a master problem and the JPS acts as a separation procedure detecting violated MaxW constraints. Instead of enforcing all MaxW constraints from the start, we progressively activate them, expecting that many constraints will be satisfied without being enforced in the model. At each iteration, the CP problem is solved considering only the currently active subset of MaxW constraints. The resulting solution gives values for variables $s_{i,j}$, $e_{i,j}$, and $d_q^k$. From this solution, we reconstruct an explicit shift-based schedule using the JPS algorithm, that schedules both processing tasks and rest intervals while respecting the makespan $c_{\max}$. Since JPS is designed to produce an optimal preemptive schedule on a single machine, we know that if the produced schedule satisfies all MaxW constraints, the solution is optimal for our whole problem. Otherwise, some deactivated MaxW constraints are violated. For each operator, we build a maximal set of non-overlapping violated intervals as follows. We repeatedly select the most violated constraint and remove from consideration all other violated constraints overlapping it, until no violated constraint remains. Selected constraints are then activated in the CP model for the next iteration. The initial pool of MaxW constraints is generated according to this same procedure, assuming that none of them are initially satisfied. This guided incremental activation significantly limits the number of MaxW constraints that are simultaneously activated.

\textbf{The baseline model} \quad For baseline comparison, we implemented a model for the MaxW-pJSP using IBM\textsuperscript{\tiny{\textregistered}} ILOG CPLEX Optimizer(CPO), a state-of-the-art CP solver for scheduling~\cite{Laborie_cpo_2018}. Since the \textsc{PreemptiveNoOverlap} constraint is not available in CPO, the model must rely on a unit-duration formulation in which each working shift of every task and each rest shift has its own decision variable. This time-indexed representation drastically increases the number of variables and enlarges the search space. The CPO model was developed following recommendations provided in the literature by one of the solver's designers~\cite{LUNARDI2020105020}. The complete model is reported in Appendix [\ref{CPO-model}].
\vspace{-0.2cm}
\section{Experiments}
\label{expes}
\vspace{-0.1cm}
To evaluate the MaxW-pJSP problem, we extended 78 pJSP benchmark instances from the literature~\cite{Juvin2023} by generating 9 MaxW instances for each pJSP instance. Each MaxW instance is generated using two parameters: (i) global density ($gd$) of rest shifts over the scheduling horizon, and (ii) local density ($ld$) of rest shifts within MaxW constraints. Randomness is introduced in both position and duration of rest intervals to avoid overly regular patterns. For each pJSP instance, a MaxW instance is generated for each density pair in the cartesian product of $\{0.1,0.25,0.4\}^2$, resulting in $78 \times 9 = 702$ MaxW-pJSP instances. Full dataset can be found on Github \cite{githubrepo} and the generation procedure is detailed in appendix [\ref{MaxW-generation}].

We allocate 45 minutes per instance. Our iterative approach, which selectively adds MaxW constraints, is currently limited by the Mistral solver that does not support dynamic constraint addition. Adding MaxW constraints requires recreating the model and restarting the search, inducing significant overhead. Despite this limitation, the method already demonstrates promising performance. We expect that a solver supporting dynamic constraint addition (e.g. TEMPO~\cite{tempo}, currently being developed) would allow for more efficient solving.

\begin{figure}[htb]
\begin{center}
    \includegraphics[scale=0.7]{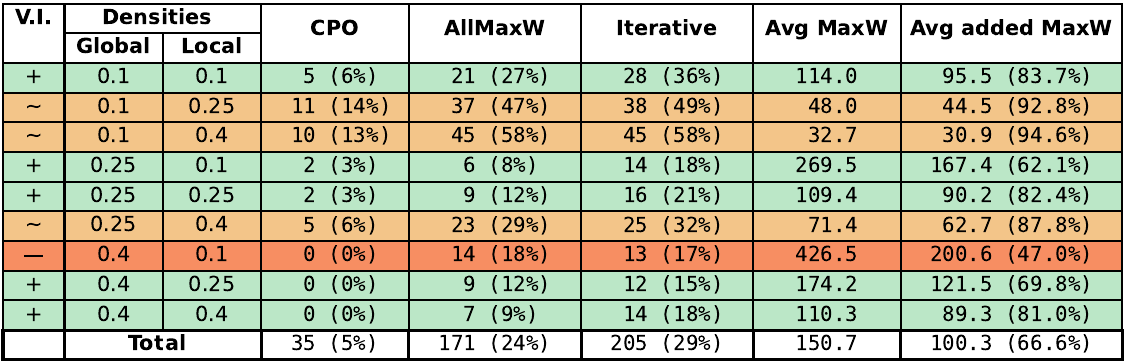}
    \captionsetup{skip=0.1cm}
    \caption{Optimal results table aggregated by density pairs (V.I. = visual impairments)}
    \label{fig:tab_density}
\end{center}
\end{figure}

Let us call \textbf{AllMaxW} the Mistral model where all MaxW constraints are activated from the beginning, \textbf{Iterative} the Mistral model with progressively activated MaxW constraints, and \textbf{CPO} the model using CPO, provided in appendix [\ref{CPO-model}]. Results unsurprisingly show that  Mistral-based methods significantly outperform CPO, solving far more instances to optimality across the entire benchmark set, as can be seen in Figure~\ref{fig:tab_density}. \textbf{Iterative} solves 205 instances to optimality (29.2\%), versus 171 (24.3\%) for \textbf{AllMaxW} and 35 (4.98\%) for \textbf{CPO} (Note that in 1 hour, only 60\% of pJSP instances are solved to optimality by the Mistral solver \cite{Juvin2023}, thus without any MaxW constraint). Our \textbf{Iterative} method performs satisfactorily: despite having to regularly restart the search from scratch, it solves 34 more instances to optimality than the full model, and is generally faster. This confirms that adding selectively MaxW constraints reduces computational effort. The poor performance of the CPO model can be attributed to the lack of a built-in \textsc{PreemptiveNoOverlap} constraint.

\vspace{-0.8cm}
Breaking down the results by density pairs (global $gd$ - local $ld$) reveals two main behaviors. Firstly, the \colorbox{my_green}{benefit of relaxation (+)}. The \textbf{Iterative} method significantly outperforms the static approach in instances with several MaxW constraints (e.g. $gd=0.4,\text{ }ld=0.25$). The last column shows that in these instances, a high proportion of them are not necessary to find a MaxW-feasible solution (e.g. 38\% can be deactivated in class $gd=0.25,\text{ }ld=0.1$). This exhibits that posting all constraints upfront may unnecessarily overload the solver. Secondly, the \colorbox{my_orange}{limitations of \textbf{Iterative} ($\sim$)}. When constraints are scarce (e.g. $gd=0.1,\text{ }ld=0.4$), both methods perform similarly: the initial \textbf{AllMaxW} model is light enough to be solved efficiently. The \colorbox{my_red}{$gd=0.4,\text{ }ld=0.1$ class (—)} highlights a bottleneck: with a massive average of 426 constraints, the \textbf{Iterative} method performs slightly worse than \textbf{AllMaxW} (1 less instance solved to optimality). Frequent violations makes the algorithm spend too much time in the iterative generation loop rather than searching. This limitation suggests that designing an adaptive constraint addition strategy is necessary for highly dense MaxW instances.

\begin{figure}[htb]
\begin{center}
    \includegraphics[width=12cm, height=5.2cm]{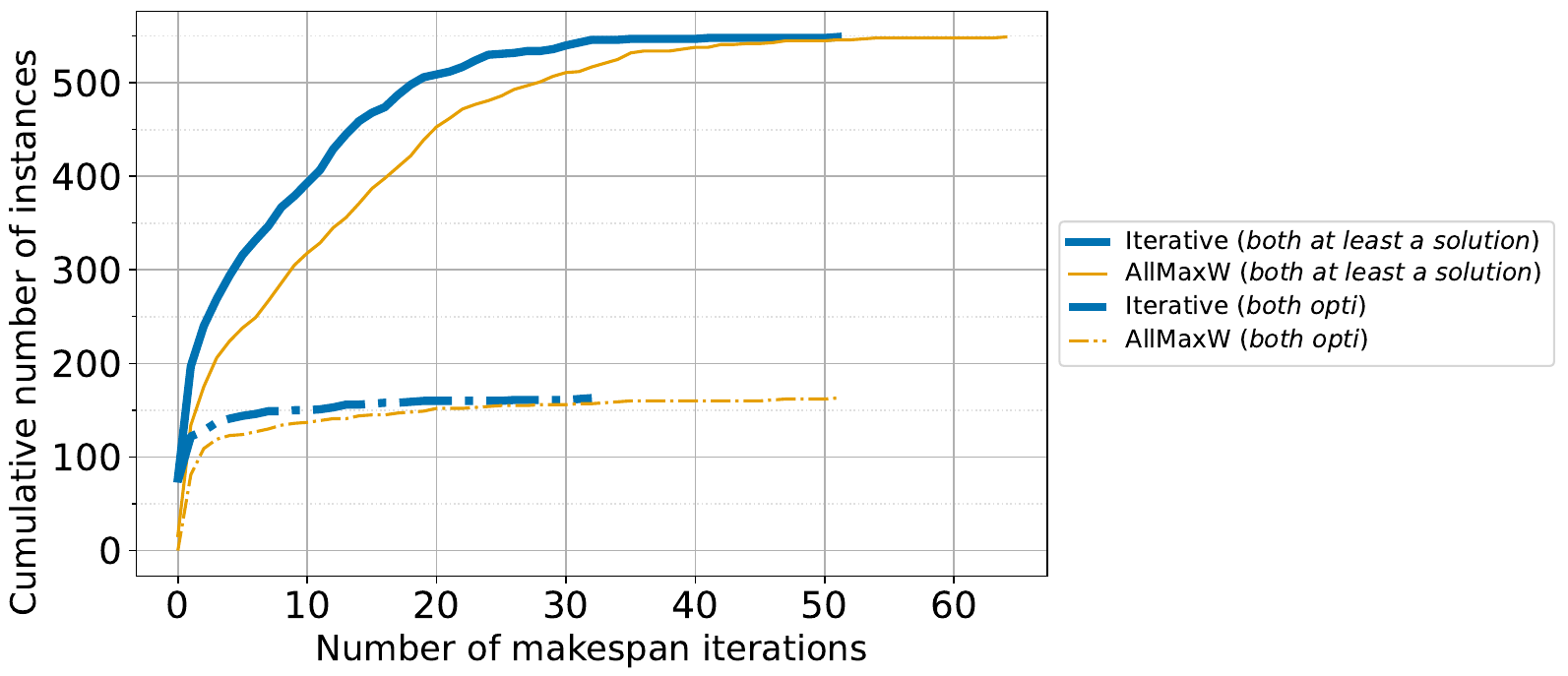}
    \captionsetup{skip=0.1cm}
    \caption{Cumulative distribution of number of makespan iterations}
    \label{fig:nb_iter}
\end{center}
\end{figure}

In addition, the number of instances solved to optimality at 20s/200s/2000s is 129/168/195 for the \textbf{Iterative} method, 113/148/167 for \textbf{AllMaxW}, and 7/12/31 for \textbf{CPO}. On the subset of instances solved optimally by both Mistral approaches, \textbf{Iterative} and \textbf{AllMaxW} perform almost identically (resp. 113/144/160 vs. 111/144/162), confirming that the \textbf{Iterative} approach does not sacrifice convergence speed on easier instances while solving strictly more hard instances overall. We reinforce these observations by showing in Figure~\ref{fig:nb_iter} that the \textbf{Iterative} method requires fewer makespan iterations (each time the solver finds a new better solution, should not be mistaken with MaxW iterations of the \textbf{Iterative} method).
\vspace{-0.2cm}
\vspace{-0.5cm}

In Figure~\ref{fig:nb_iter}, solid curves only consider instances where \textbf{AllMaxW} and \textbf{Iterative} methods found a solution, and dashdot curves only consider instances where both methods proved optimality. Indeed, considering all instances would be unfair to \textbf{AllMaxW} that finds way less feasible solutions and optimality proofs. For example, Figure~\ref{fig:nb_iter} shows that with 15 makespan iterations, our \textbf{Iterative} approach finds a solution for 468 instances (66.7\%), and proves optimality for 157 of them (22.3\%). This outperforms the \textbf{AllMaxW} method that finds a solution for 387 instances (55.1\%), amongst which 145 (20.6\%) are shown optimal. This highlights the efficiency and practical relevance of the \textbf{Iterative} approach.

\begin{figure}[htb]
\begin{center}
    \includegraphics[width=12cm, height=6cm]{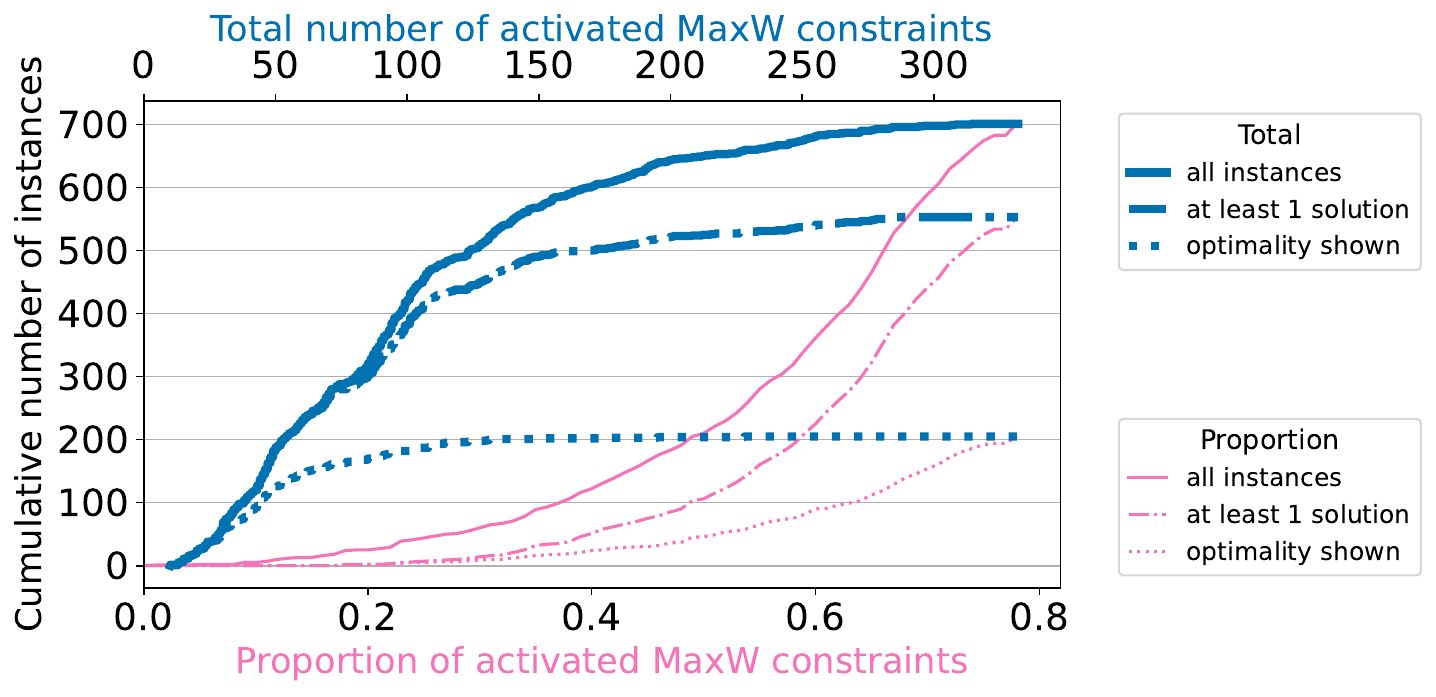}
    \caption{Cumulative distribution of activated MaxW constraints}
    \label{fig:prop_MaxW}
\end{center}
\end{figure}

To better understand the performance gap between the two methods, we investigate the actual benefit of saving MaxW constraints using the iterative approach. Figure~\ref{fig:prop_MaxW} reports the \colorbox{my_pink}{proportion (thin)} and \colorbox{my_blue}{total count (thick)} of activated MaxW constraints for: (i) all instances, (ii) instances with at least one feasible solution found, (iii) instances where optimality was proven.

First, the \colorbox{my_pink}{proportion (thin)} curves clearly demonstrate that a vast majority of the constraints are never activated during the search. For example, in more than 500 instances, the solver finds a feasible solution with less than 70\% of the MaxW constraints activated (dashdot thin pink curve). However, looking at the dotted thin pink curve, one can note that our method struggle to prove optimality without a high proportion of MaxW added.

Furthermore, the \colorbox{my_blue}{total count (thick)} curves reveal a significant plateau effect: once more than $\sim$100 MaxW constraints are activated, the number of instances for which a feasible solution is found strongly stagnates. This plateau is even more pronounced when it comes to proving optimality. This observation strongly justifies the core principle of our \textbf{Iterative} method, demonstrating that injecting too many constraints simultaneously overburdens the solver, suggesting the need to precisely target a subset of necessary and sufficient MaxW.
\newpage
\vspace{-0.5cm}
These results suggest an important avenue for future research: comparing different MaxW-selection heuristics to work towards a better incremental activation of MaxW constraints. We indeed showed with previous figures that it is of utmost importance for providing high-quality relaxations, when a lot of MaxW constraints highly increases complexity. This iterative \textit{lazy constraint generation} could be extended to multiple problems that are very hard to solve. It is important to remind that our \textbf{Iterative} approach is quite limited by the Mistral software, that requires restarting the search everytime that a MaxW constraint is added. It is still a very promising approach, giving best results than the \textbf{AllMaxW} classical method.

\vspace{-0.2cm}
\section{Conclusion and future work}
\vspace{-0.1cm}
In this paper, we proposed and validated a lazy constraint generation approach with a reconstruction-based separation procedure to efficiently handle maximum workload constraints in pJSP. In addition, we introduced a new benchmark composed of 702 MaxW-pJSP instances derived from classical pJSP benchmarks, covering a wide range of global/local rest-shift densities. Our experiments show that solver performance is strongly influenced by the number of MaxW constraints. Indeed, benefits of the iterative approach become pronounced as the number of constraints increases and as the density of required rest shifts grows. In workforce scheduling practice, such constraints often arise in large quantities due to rolling rules (e.g. “no more than five working shifts in any seven-day window”). In such contexts, our approach could clearly outperform static formulations, which activates all constraints from the start.

Several research directions emerge from this work. First, the implementation is limited by the Mistral solver, which requires frequent search restarts and model recreation when constraints are added, thereby introducing significant overhead. Future work will therefore explore the use of the TEMPO solver~\cite{tempo}, where the \textsc{PreemptiveNoOverlap} constraint is currently being integrated. Since TEMPO natively supports dynamic constraint addition, it should enable a significantly more efficient implementation of the iterative method. Finally, we aim to extend our framework to incorporate minimum workload (MinW) constraints, ensuring that employees are assigned a sufficient amount of work. Preliminary investigations show that MinW cannot be modeled as a simple symmetric counterpart of MaxW, since enforcing a maximum number of rest shifts does not directly translate into enforcing a minimum number of working shifts, due to possible idle shifts that are not controlled by the model but by the JPS algorithm. Designing elegant models and efficient solution methods for this extension constitutes an important direction for future research.

\newpage
\bibliography{lipics-v2021-sample-article}

\appendix

\section{Appendix A - CPO model}
\label{CPO-model}

In CPO, without preemption, modeling is simple and efficient, as tasks are single fixed intervals and strong filtering algorithms (e.g., edge-finding~\cite{Vilim_filtering_2005}) apply directly. With preemption, however, each task must be split into multiple unit-length optional intervals, dramatically increasing combinatorial complexity and making propagation more difficult.

The extensive CPO model used for experiments is the following (note that when an optional interval is absent, it is not belonging to any interval). Explanations follow the model. In the following, if not specified, consider $i\in \{0, ..., |\mathcal{J}|-1\}$, $j\in \{0, ..., n_i\}$, $k \in \textsc{Operators}$. We write $\delta_c, u_c,v_c$ the triplet values of MaxW constraint $c$.

$$\min \quad \max_{i,j} \big(\operatorname{endOf}(W_{i,j,P_{ij}})\big)$$
\begin{align*}
\text { interval } W_{ij\ell},\, \operatorname{size}=1
&& \quad \forall i, \forall j, \forall \ell \in \{1,\text{ }, P_{ij}\}\quad (\text{V}1) \\
\text { interval } O_{kc\lambda},\text{ optional },\operatorname{size}=1
&& \quad \forall k, \forall c \in \text{MaxW}_k, \forall \lambda \in \{1,\text{ }, \delta_c\}\quad (\text{V}2) \\
O_{kc\lambda} \subseteq [u_c, v_c]
&& \quad \forall k, \forall c \in \text{MaxW}_k, \forall \lambda \in \{1, \text{ }, \delta_c\}\quad (\text{C}1) \\
\textsc{EndBeforeStart}(W_{ij,\ell-1}, W_{ij,\ell})&&\quad \forall i, \forall j, \forall \ell \in \{2,\text{ }, P_{ij}\}\quad (\text{C}2) \\
\textsc{EndBeforeStart}\left(W_{i,j,p_{ij}}, W_{i,j+1,0}\right)&& \quad \forall i, \forall j \in \{0, \text{ } , n_i -1\}\quad (\text{C}3)  \\
\textsc{NoOverlap}\left(\underset{\forall i, j \in \text{T}_k, \ell}{\left[W_{ij\ell}\right]} \cup \underset{\forall c\in \mathit{MaxW}_k, \lambda}{\left[O_{kc\lambda}\right]}\right) &&\quad \forall k \quad (\text{C}4) \\
\delta_c \leq \sum_{\lambda=1, ..., \delta_c}\big( O_{kc\lambda} \in [u_c, v_c]\big) && \quad \forall k, \forall c \in \text{MaxW}_k \quad (\text{C}5) 
\end{align*}

In this model, we have two \textit{families} of decision variables. For each task $j$ of each activity $i$, we add as many unit-length interval variables (V1) as the processing time of the task (indexed by $\ell$). To deal with MaxW constraints, the model also features a set of interval variables $O$ (for rest shifts). For every MaxW constraint $c$, we add to operator $k$ $\delta_c$ unit-length \textit{optional} interval variables (V2). If these optional variables are included in the final solution, they must be scheduled within the interval $[u,v]$ (C1). These variables allow flexibility in how rest periods are scheduled to satisfy workload regulations. The objective is still to minimize the makespan. Constraint C2 enforces sequential execution of a task: for each unit $\ell$ of any task $ij$ except the last, the shift represented by variable $W_{ij\ell}$ must end before $W_{ij(\ell+1)}$ begins. Precedence constraints between successive operations of the same activity are also imposed (C3). Constraint C4 enforces a no-overlap condition, ensuring that for each operator, all its working ($W$) and rest ($O$) shifts do not overlap in time. Lastly, this model handles the maximum workload requirements (C5). For each operator $k$ and each constraint MaxW$_{k}$ $c$ defined over a time window $[u, v]$, the model requires for at least $\delta$ optional intervals to be present in this window.

Note that we tried adding symmetry-breaking constraints for the rest shifts variables (intervals $O_{kc\lambda}$), and that gave worst results. Indeed, their ordering makes the search space smaller, but also more constrained, which is de facto not beneficial for our model.

\section{Appendix B - MaxW instances generation}
\label{MaxW-generation}

The first main parameter of our algorithm to generate MaxW instances is the global desired off shift density $\mathcal{D}$. It expresses the minimal overall proportions of rest days that each operator should have over the entire horizon ($\mathcal{H}$). These rest days are distributed across multiple time intervals, which are randomly chosen but must respect some constraints. Each interval should be at least 3 days, and at most $0.2*\mathcal{H}$ days (neither too short nor too long). Within each interval, the number of rest days is related to the other main parameter, the local minimal density ($d$). They express the target proportion of rest days that the operator should have over any specific MaxW interval. These proportions are chosen with random gaussian probability: $\mathcal{N}(d, d/4)$. This method introduces some unpredictability into when rest periods occur and how long they are. MaxW constraints are incrementally generated for each operator until the total minimal number of rest days for all operators meet the target.
\end{document}